\documentclass{siamltex}

\usepackage{multirow}
\usepackage{caption}
\usepackage{latexsym,graphicx, amssymb, amsfonts, subfigure, setspace, geometry, amsmath}
\usepackage{pgfplots}
\usepackage{grffile}
\usepackage{pgfplotstable}
\usepackage{mathtools}
\usepackage[utf8]{inputenc}
\usepackage[english]{babel}

\pgfplotsset{compat=newest}

\newtheorem{remark}{Remark}[section]

\newcommand{\xb}{\mathbf{x}}
\newcommand{\yb}{\mathbf{y}}

\newcommand{\db}{\mathbf{d}}

\newcommand{\parB}{{\partial \mathcal{B}}}

\title{High resolution inverse scattering in two dimensions using recursive
linearization}

\author{Carlos Borges\thanks{Courant Institute of Mathematical Sciences, New York University, New York, NY}
\and
Adrianna Gillman\thanks{Computational and Applied Mathematics, Rice University, Houston, TX}
\and 
Leslie Greengard\thanks{Courant Institute of Mathematical Sciences, New York University, New York, NY and Simons Center for Data Analysis, Simons Foundation, New York, NY}
}

\begin{document}

\maketitle

\begin{abstract}
We describe a fast, stable algorithm for the solution of the inverse acoustic
scattering problem in two dimensions. Given full aperture
far field measurements of the scattered field for multiple angles of 
incidence, we use Chen's method of recursive linearization
to reconstruct an unknown sound speed at resolutions of 
thousands of square wavelengths in a fully nonlinear regime.
Despite the fact that the underlying optimization problem is formally 
ill-posed and non-convex, recursive linearization requires only the
solution of a sequence of linear least squares
problems at successively higher frequencies. 
By seeking a suitably band-limited approximation of the sound speed profile, 
each least squares calculation is well-conditioned and
involves the solution of a large number of forward scattering problems, 
for which we employ a recently developed, spectrally accurate, 
fast direct solver. For the largest problems considered, involving 19,600
unknowns, approximately one million partial differential equations were solved,
requiring approximately two days to compute using a parallel MATLAB implementation
on a multi-core workstation.

\end{abstract}

\section{Introduction}
Inverse scattering problems arise in many areas of science and engineering, 
including medical imaging \cite{kuchment2014radon,
2006-MMRAMI-WM,nashed2002inverse,scherzer2010handbook}, remote sensing 
\cite{Ustinov2014,Wang01022012}, ocean acoustics \cite{chavent2012inverse,0266-5611-10-5-003}, 
nondestructive testing \cite{collins1995nondestructive,engl2012inverse}, 
geophysics \cite{aster2013parameter,tarantola2013inverse} and radar 
\cite{BorgesG15,cheney2009fundamentals,Colton}. 
In this paper, we investigate the problem of recovering an unknown
compactly supported sound speed profile or {\em contrast function}, 
denoted by $q(\xb)$, from far-field acoustic scattering measurements in 
two space dimensions.

Letting $\Omega$ denote a domain containing the support of $q(\xb)$, 
we very briefly review the
forward scattering problem in the time-harmonic setting, when
the contrast function is known. The governing equation is then the 
Helmholtz equation  
\begin{equation}
  \Delta u(\xb) + k^2 (1 - q(\xb)) \, u(\xb) = 0, 
\label{origpde}
\end{equation}
for $\xb \in \mathbb{R}^2$, where 
\[ u(\xb) = u^{\emph{inc}}(\xb) + u^{\emph{scat}}(\xb)
\]
and $k$ is the frequency (or wavenumber) under consideration.
Here, 
$u^{\emph{inc}}$ denotes a known incoming field, which satisfies the constant
coefficient Helmholtz equation 
\begin{equation}
\Delta u(\xb) + k^2 \, u(\xb) = 0, 
\label{hhom}
\end{equation}
and $u^{\emph{scat}}$ denotes the unknown scattered field, which 
must satisfy the Sommerfeld radiation condition
\begin{equation}
\lim_{r \rightarrow \infty}
\sqrt{r} \, \left( \frac{\partial u^{\emph{scat}}}{\partial r} -
iku^{\emph{scat}} \right) = 0,
\label{srad}
\end{equation}
where $r=\|\xb\|$.
It is straightforward to verify that
\begin{equation}
\Delta u^{\text{\emph{scat}}}(\xb) + k^2 (1-q(\xb)) u^{\text{\emph{scat}}}(\xb) 
= k^2q(\xb) u^{\text{\emph{inc}}}(\xb) \, ,
\label{fscatprob}
\end{equation}
which reduces to the constant coefficient equation
\eqref{hhom} outside the support of $q(\xb)$.
Together, \eqref{fscatprob} and \eqref{srad} define the 
forward scattering problem.

\begin{figure}
\center
\includegraphics[width=2in]{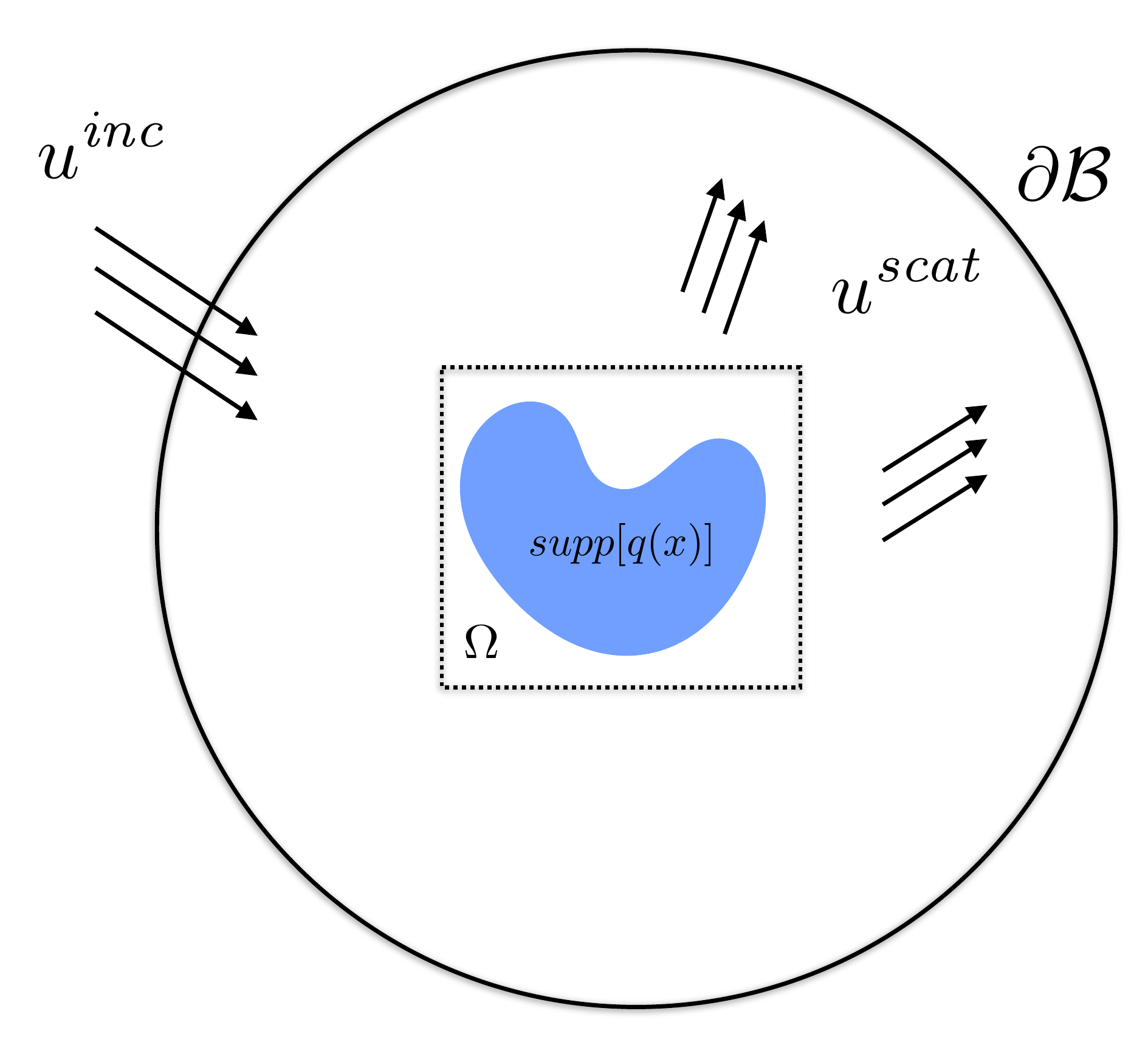}
\caption{Scattering from a compact inhomogeneity 
(contrast function) $q(\xb)$: the support of $q(\xb)$ 
is assumed to lie within
the domain $\Omega$ and impinged upon by 
an incoming field $u^{\emph{inc}}$, such as
a plane wave. 
In the {\em forward
scattering problem}, $q(\xb)$ is known and 
one seeks to compute the scattered field, 
either within
$\Omega$ or in the far field - say, on the boundary $\partial\mathcal{B}$ 
of an enclosing disk $\mathcal{B}$. 
In the {\em inverse scattering problem}, $q(\xb)$ is unknown, and one seeks
to determine it from measurements of the scattered field on
$\partial\mathcal{B}$.
}\label{fig1}
\end{figure}

We assume that the incoming field is a plane wave of the form
\[ u^{\emph{inc}}(\xb)=\exp(i k \, \xb \cdot \db), \]
where $\db$ is a unit vector that defines the direction of propagation.
We also assume that the scattered field is measured on the boundary 
$\partial \mathcal{B}$ of a disk $\mathcal{B}$ which contains $\Omega$ 
(Fig. \ref{fig1}). More precisely, we denote by 
$u^{far}(\theta)$ the measured data
\[ u^{far}(\theta) = 
u^{scat}(R \cos \theta, R \sin \theta), \]
for $\theta \in [0,2\pi]$, where $R$ denotes the radius of the disk $\mathcal{B}$.

\vspace{.1in}

\begin{remark}
When it is important to be explicit about the direction of incidence and frequency,
we will denote $u^{\emph{inc}}(\xb)$ by 
$u^{\text{\emph{inc}}}_{k,\db}(\xb)$ and 
$u^{\emph{scat}}(\xb)$ by 
$u^{\text{\emph{scat}}}_{k,\db}(\xb)$.
Likewise, when it is necessary to be explicit about the dependence on $q(\xb)$,
we will denote $u^{\emph{scat}}(\xb)$ by 
$u_q^{\emph{scat}}(\xb)$ and 
$u^{\text{\emph{scat}}}_{k,\db}(\xb)$ by
$u^{\text{\emph{scat}}}_{q,k,\db}(\xb)$.
The scattered field measured on $\parB$ will be denoted by
$u^{far}(\theta)$ or
$u^{far}_{k,\db}(\theta)$.
\end{remark}

\vspace{.1in}

\vspace{.1in}

\begin{definition}
Suppose that, for a fixed frequency $k$, a series of experiments
is carried out, with $M$ distinct
plane waves impinging on a domain $\Omega$ which contains
the support of an unknown contrast function $q(\xb)$.
Let the incident directions be denoted by 
$\{ \db_{m}, m=1,\dots,M \}$.
The {\em single frequency inverse scattering problem} consists of 
determining $q(\xb)$ from 
$\{ u^{{far}}_{k,\db_{m}}(\theta), m = 1,\dots,M \}$.
\end{definition}

\vspace{.1in}

It is important to note that, in the far field, no more than
$O(k)$ independent measurements can reasonably be made on 
$\partial \mathcal{B}$,
assuming the support of $q(\xb)$ has been normalized to have approximately
unit diameter. This follows
from either standard estimates for the behavior of 
the multipole expansion of $u^{\emph{scat}}_{k,\db_{m}}(\xb)$, or
the Heisenberg uncertainty principle 
\cite{Chenrep1088,Chenrep1091,Chen,crutchfield2006remarks}.
In physical terms, the issue is that Fourier modes on $\partial \mathcal{B}$
whose frequency exceeds $k$ correspond to evanescent and rapidly
decaying fields emanating from the scatterer. Acquiring such data
would impose exponential accuracy requirements on the measurements of 
$u^{far}_{k,\db_{m}}(\theta)$.
In short,
only $O(k)$ linearly independent measurements are available
for each angle of incidence with finite precision. 
Similar arguments show that only $O(k)$ 
independent directions of incidence are useful in probing the unknown 
inhomogeneity, leading to a total of $O(k^2)$ independent measurements. 
Thus, in two dimensions, the single frequency inverse problem is at the 
limits of feasibility 
in seeking to reconstruct a model for $q(\xb)$ with $O(k^2)$
unknowns.

\vspace{.1in}

\begin{definition}
Suppose now that we probe the unknown function $q(\xb)$ at a set
of frequencies $\{ k_j, \, j = 1,\dots,Q\}$, with
incident directions at each frequency $k_j$ denoted by 
$\{ \db_{j,m}, \, m=1,\dots,M_j \}$.
The {\em multi-frequency inverse scattering problem} consists of determining 
$q(\xb)$ from 
$\{ u^{far}_{k_j,\db_{j,m}}(\theta);\, j = 1,\dots,Q,\, 
m = 1,\dots,M_j \} $.
\end{definition}

\vspace{.1in}

\begin{remark}
As indicated above, the number of linear independent measurements that can be 
made on $\parB$ is of the order $O(k_j)$ at frequency $k_j$. We will
denote by $P_j$ the number of distinct (equispaced) measurements made in the 
angular variable $\theta$.
In practice, one could make a larger number of measurements and filter/denoise
the data by fitting a Fourier series on $\parB$ with $P_j$ modes.
\end{remark}

\vspace{.1in}

We assume $q(\xb) \in C_{0}(\Omega)$ and
define the operator 
$\mathcal{F}_{k,\db}:C_{0}(\Omega)\rightarrow L^2(\partial \mathcal{B})$ 
by
\begin{equation}
  \mathcal{F}_{k,\db} [q] = u^{\emph{far}}_{k,\db}.
\label{farfielddef}
\end{equation}
The operator $\mathcal{F}$ is well-defined since the forward scattering
problem is well-posed.
To obtain the value of $u^{far}_{k,\db}$ at a point 
$\xb = (R \cos \theta, R \sin \theta)$, 
one must solve 
\eqref{fscatprob} and \eqref{srad} or
its integral equation counterpart, the Lippmann-Schwinger equation
\cite{Colton,Nedelec},
\begin{equation}
u_{q,k,\db}^{\text{\emph{scat}}}(\xb)+k^2\iint_\Omega G(\xb,\yb) q(\yb) (u_{q,k,\db}^{\text{\emph{scat}}}(\yb)+u^{\text{\emph{inc}}}(\yb,k,\db))~d{\yb}=0, 
\label{inteqscat}
\end{equation}
with $G(\xb,\yb) = \frac{i}{4} H_0^{(1)}(k \|\xb-\yb\|)$ where 
$H_0^{(1)}(x)$ is the usual Hankel function of the first kind.
Eq. \eqref{inteqscat} is derived by integrating both sides of
\eqref{fscatprob} against $G(\xb,\yb)$, 
using the fact that it is the Green's function for 
eq. \eqref{hhom} satisfying the radiation condition \eqref{srad}.

We will focus here on the multi-frequency inverse scattering 
problem defined above. In other words, our goal is to solve the nonlinear
system of equations
\begin{equation}
 \mathcal{F}_{k_j,\db_m} [q] = 
u^{\emph{far}}_{k_j,\db_{j,m}},
\label{invscatmf}
\end{equation}
for $j=1,\dots,Q$, 
$m=1,\dots,M_j$. 
This is an ill-posed, nonlinear and nonconvex problem with a substantial
literature (see, for example, \cite{BaoLi2015,Colton,kress2007} and the references therein).
Broadly speaking, existing approaches can be classified as 
either iterative methods, derived from a nonlinear optimization framework,
or direct methods, based on ideas drawn from image and signal processing. 
Iterative methods include
variants of Newton's method 
\cite{Chenrep1088,Chenrep1091,Chen}, the Gauss-Newton method 
\cite{buithanh2012analysisobj,buithanh2012analysis,buithanh2013analysis}, 
Landweber iteration \cite{Bao_li_2007,Bao_li_2005_1,Bao_li_2005,Bao_li_2009,li_bao,Bao,Hohage2001}, quasi-Newton methods \cite{GUTMAN19935, 
GutmanKlibanov93,GutmanKlibanov94}, and the nonlinear 
conjugate gradient method \cite{KLEINMAN199217,kleinmanvdb93,vandenberg97}.
Direct methods include decomposition methods 
\cite{ColtonKirsch91,colton1988inverse,kirsch_potential,
Potthast_apoint-source,Potthast_apoint-source1,potthast2001point}, 
the linear sampling method \cite{CakoniColtonMonk,Colton_kirsch}, 
the singular source method \cite{Potthast_single,potthast2001point},
the factorization method \cite{kirsch1996introduction,Kirsch_2}, 
and the probe method of Ikehata \cite{Ikehata}. 
Nevertheless, most numerical work on reconstruction has been limited
to fairly simple contrast functions involving perhaps dozens 
of parameters in a model for the unknown contrast function $q(\xb)$.

In this paper, we are interested in developing a method for high-resolution
two-dimensional applications, where $q(\xb)$ is modeled as a
function on a grid with up to 
$100 \times 100$ unknowns. We will make use 
of a Newton-like iterative method which relies on the frequency $k$ 
as a continuation parameter.
More precisely, we will solve a sequence of {\em single-frequency} inverse
problems for higher and higher values of $k$,
using the approximation
of $q(\xb)$ obtained at the preceding frequency as an initial guess.
In the context of inverse scattering, such a scheme was first proposed 
by Chen \cite{Chen} and is referred to as {\em recursive linearization}.
More recent contributions include
\cite{Bao_li_2005_1,Bao_li_2005,Bao_li_2009,Bao}.
The analogous problem for scattering from an unknown, impenetrable,
sound-soft object is discussed in
\cite{BorgesG15,Sini2013,Sini11}.
For time-domain versions of the problem, see
\cite{Beilina2015,beilina2015-1}.

\vspace{.1in}

\begin{remark}
The ill-posedness inherent in inverse scattering is closely tied to the
issues stemming from the Heisenberg uncertainty principle discussed above.
Loosely speaking, features of $q(\xb)$ that have frequency content greater
than the probing incident field are evanescent and poorly determined by 
far field measurements. Overcoming this problem is 
often addressed by using some form of {\em ad hoc} regularization while 
solving the linearized subproblems which arise in the various reconstruction
schemes \cite{cakoni2006,Colton,kaipio2010,kirsch_book}. 
In the original work on recursive linearization
\cite{Chenrep1088,Chenrep1091,Chen}, however, and in our previous work on 
inverse obstacle scattering \cite{BorgesG15}, it was shown that the same 
stabilizing effect can be achieved by using a suitably band-limited
model for the unknown. We will continue to employ that strategy here
(see Section \ref{sec:inv}).
\end{remark}

\vspace{.1in}


An outline of the paper follows. In Section \ref{sec:direct},
we describe the forward scattering problem and its solution using the 
fast direct solver developed in \cite{Gillman2014} - the so-called
Hierarchical Poincar\'{e}-Steklov method.
In Section \ref{sec:inv}, we describe our implementation of recursive
linearization for the inverse problem and
in Section \ref{sec:numerics}, we illustrate the performance of our
method. Section \ref{sec:concl} contains some concluding remarks
and a discussion of future directions for research.

\section{The direct scattering problem}
\label{sec:direct}

In this section, we briefly review the forward scattering problem
and its solution for penetrable media in two dimensions. 
We assume that the index of refraction 
$1-q(\xb)$ is real and positive for $\xb \in \Omega$, so that the problem 
has a unique solution for any $k>0$ \cite{Colton}.

We begin by observing that an alternative formulation for the 
original partial differential equation 
\eqref{origpde} is to consider an {\em interior} variable
medium problem
\begin{eqnarray}
\Delta u(\xb ) +k^2(1-q(\xb)) u(\xb ) &=&0 \quad  \mbox{in} \quad \Omega, \nonumber \\
u(\xb )&=&h(\xb ) \quad  \mbox{on} \quad \partial\Omega,
\label{intdir}
\end{eqnarray}
coupled with an {\em exterior} constant-coefficient problem
\begin{eqnarray}
\Delta u^{scat}(\xb ) +k^2 u^{scat}(\xb ) &=&0 \quad \quad \mbox{in} \quad \mathbb{R}^2\setminus\Omega,\\
u^{scat}(\xb )&=&s(\xb ) \quad \quad \mbox{in} \quad \partial\Omega,\nonumber \\
 \frac{\partial u^{scat}}{\partial r} -iku^{scat} &=& o(r^{-1/2})\quad  r=\|\xb \|\rightarrow\infty.
\label{extdir}
\end{eqnarray}
For the sake of simplicity, let us assume
that the interior Dirichlet problem does not have a resonance at the 
particular frequency $k$ under consideration. 
We then seek to find functions $h(\xb)$ and $s(\xb)$ so that 
gluing together the interior and exterior total fields yields a
continuously differentiable {\em total field} $u(\xb)$. 
If that can be achieved, 
then the solution to \eqref{intdir} matches the solution to
\eqref{origpde} in the interior of $\Omega$
and $u = u^{scat} + u^{inc}$ matches the solution to
\eqref{origpde} in the exterior of $\Omega$ by a simple uniqueness argument 
\cite{Colton}.

To accomplish this matching, let $\frac{\partial u}{\partial n}$ denote
the outward normal derivative of the solution to \eqref{intdir} on
$\partial\Omega$. We may then define the interior
``Dirichlet-to-Neumann" map $T^{int}$ by
\[ T^{int} h = \frac{\partial u}{\partial n}. \]
There is also a well-defined
exterior ``Dirichlet-to-Neumann" map $T^{ext}$ such that 
\[ T^{ext} s = \frac{\partial u^{scat}}{\partial n}. \]
Given these two maps, it is straightforward to determine $s(\xb)$ and
$h(\xb)$ by impose the continuity conditions
\begin{align*}
s(\xb) + u^{inc}(\xb) &= h \\
T^{ext}[s](\xb) + \frac{\partial u^{inc}}{\partial n}(\xb) &= T^{int}[h](\xb).
\end{align*}
In particular , we can obtain the scattered field $u^{scat}(\xb) = s(\xb)$ 
on $\partial\Omega$ by solving the problem 
(analogous to equation (2.12) in\cite{kirschmonk94}):
\begin{equation*}
\left(T^{int}-T^{ext}\right) u^{scat}\vert_{\partial\Omega} = \frac{\partial u^{inc}}{\partial n} - T^{int} u^{inc}\vert_{\partial\Omega}.
\end{equation*}

\vspace{.1in}

\begin{remark}
While $T^{int}$ is rather complicated to describe, 
$T^{ext}$ can be written using standard layer potentials from
Green's formula, since the scattered field
$u^{scat}(\xb)$ satisfies
\begin{equation*}
u^{scat}(\xb)=Du^{scat}(\xb)-S \frac{\partial u^{scat}}{\partial n} (\xb)
\end{equation*}
for $\xb$ in the exterior of $\Omega$.
Here, $D\phi(\xb)=\int_{\partial\Omega}\frac{\partial G(\xb,\yb)}{\partial n_{y}}
\phi(\yb)ds(\yb)$ and $S\phi(\xb)=\int_{\partial\Omega}G(\xb,\yb)\phi(\yb)ds(\yb)$
are the double and single layer operators, respectively and 
$G(\xb,\yb)=(i/4)H_0^{(1)}(k\|\xb-\yb\|)$. Using standard jump relations
\cite{colton1983integral,Colton}, it is easy to verify that 
\begin{equation*}
T^{ext}=S^{-1} \left(D-\frac{I}{2} \right).
\end{equation*}
\end{remark}

\vspace{.1in}

This is the essence of the approach used in the 
Hierarchical Poincar\'{e}-Steklov (HPS) solver of \cite{Gillman2014}. 
Without entering into details, we simply note here that 
the basic discretization, for $K$th order accuracy,
involves superimposing a quad-tree on 
the domain $\Omega$, with 
tensor product $K \times K$ Chebyshev grids on each leaf node, used to
represent both $u(\xb)$ and $q(\xb)$. 
The HPS method solves the interior problem on
$\Omega$ by a recursive merging procedure, and represents
the exterior field using a layer potential on $\partial\Omega$.
By a careful use of ``impedance-to-impedance" maps, the method involves
well-conditioned operators and requires only $O(N^{3/2})$ work for factoring
the system matrix with a given contrast function $q(\xb)$. 
Given that factorization,
the solver requires only $O(N\log N)$ work in order to solve
eq. \eqref{fscatprob} for each right-hand side defined by $u^{inc}$.
(See the original paper \cite{Gillman2014} for a complete description of the 
method.)

Over the last decade, a number of fast direct solvers
have been developed with the same basic complexity, some using 
direct discretization of the partial differential equation (PDE) and some
using the Lippmann-Schwinger integral formulation.
We will not attempt to review the literature here and 
refer the reader to \cite{BorgesLS,ambikasaran2013fastdirect,
borm2003hierarchical,borm2003introduction,
chandrasekaran2006fast1,chen2002fast,corona2015,ifmm_darve,
hackbusch2001introduction,ho2012fast,
martinsson2013,xia2010fast,
leonardo2016fast} 
and the references therein.

The HPS solver was implemented in MATLAB and run in parallel mode
using up to 12 cores of a system with
2.5GHz Intel Xeon CPUs.
To illustrate its performance, Table \ref{forward_tab}
presents the run-time for a sequence of problems with increasing 
frequency $k$ and an increasing
number of discretization points, using the simple contrast function 
\begin{equation*}
q(\xb)= q(x,y)=1.5\exp\left(-\frac{x^2+y^2}{50}\right)
\end{equation*}
in the domain $\Omega = [-\frac{\pi}{2},\frac{\pi}{2}]^2$.
$N$ denotes the total number of points used to discretize the 
domain $\Omega$, and $N_{\partial\Omega}$ is the number of points used 
on the boundary $\partial \Omega$ for the solution of the exterior problem.
$T_\text{interior}$, $T_\text{bdry}$ and $T_\text{solve}$ are the times 
(in seconds) to factor the interior system matrix, the exterior system
matrix, and apply the resulting inverse, respectively.

\begin{table}
\begin{center}
\begin{tabular}{c|c|c|c|c|c}
\hline
\hspace{0.2cm}$k$\hspace{0.2cm} &
\hspace{0.2cm}$N$\hspace{0.2cm} &
\hspace{0.2cm}$N_{\partial\Omega}$\hspace{0.2cm} &
\hspace{0.2cm}$T_\text{interior}$\hspace{0.2cm} &
\hspace{0.2cm}$T_\text{bdry}$\hspace{0.2cm} &
\hspace{0.2cm}$T_\text{solve}$\hspace{0.2cm} \\
\hline\hline
1    &   3721       &  640   &          1.79e+00  &       1.94e+00    &      9.56e-04 \\
2    &   3721       &  640   &          8.79e-01   &       1.70e+00    &      1.24e-03 \\
4    &   3721       &  640   &          8.57e-01   &       1.71e+00    &      1.74e-03 \\
8    &   14641     &  800   &          4.12e+00  &       2.24e+00    &      1.47e-03 \\
16  &   58081     &  1120  &         1.66e+01  &       3.38e+00    &      2.95e-03 \\
32  &   231361   &  1760  &         6.43e+01  &       5.96e+00    &      8.70e-03 \\
64  &   923521   &  3040  &         2.66e+02  &       1.23e+01    &      2.16e-02 \\
128&   3690241 &  5600  &         1.10e+03  &       3.56e+01    &      8.71e-02 \\
\hline
\end{tabular}
\end{center}
\caption{Run times for the direct scattering problem with 
16 points per wavelength.}\label{forward_tab}
\end{table}

\begin{remark}
 Note that the performance of the solver is independent of the wavenumber.  Here the 
 number of points per wavelength is kept fixed for consistency with experiments later 
where this choice guarantees a specific accuracy.
\end{remark}

\section{The inverse scattering problem} \label{sec:inv}

We turn now to the problem of recovering $q(\xb)$ from a set of far-field
measurements of the scattered field.
Instead of solving the full multi-frequency system of equations \eqref{invscatmf}, 
we will proceed by solving a sequence of single frequency inverse problems.
At each fixed frequency $k$, we assemble the scattered data for 
each of $M$ incident directions into the nonlinear system:
\begin{equation}
{\bf F}_{k} [q] = {\bf u}_k^{far},
\label{singfreqsys}
\end{equation}
where
\begin{equation}
\arraycolsep=1.2pt\def\arraystretch{1.4}
 {\bf F}_{k} [q] \equiv
\left[ \begin{array}{c}  
\mathcal{F}_{k,\db_{1}}[q] \\
\mathcal{F}_{k,\db_{2}}[q] \\
\dots \\
\mathcal{F}_{k,\db_{M}}[q]
\end{array} \right] ,\quad
{\bf u}_k^{far}(\theta) \equiv
\left[ \begin{array}{c}  
u^{\emph{far}}_{k,\db_{1}}(\theta) \\
u^{\emph{far}}_{k,\db_{2}}(\theta) \\
\dots \\
u^{\emph{far}}_{k,\db_{M}}(\theta)
\end{array} \right].
\end{equation}

\subsection{Linearization} \label{sec:pert}

Using Newton's method, we
linearize the problem \eqref{singfreqsys} for 
$q(\xb)$ in the neighborhood of an initial guess $q_0(\xb)$. For this, let
$\delta q = q-q_0$, so that we may write
\begin{equation}
{\bf F}_{k}[q_0]+{\bf J}_{q_0,k} \, \delta q \approx 
{\bf F}_{k}[q_0 + \delta q] = {\bf u}_k^{far},
\end{equation}
leading to the linear system
\begin{equation}
{\bf J}_{q_0,k} \, \delta q =
{\bf u}_k^{far}
- {\bf F}_{k}[q_0],
\label{linsingfreqsys}
\end{equation}
where ${\bf J}_{q_0,k}$ is the Fr\'{e}chet derivative of the 
operator ${\bf F}$ at $q_0$: 
\begin{equation}
\arraycolsep=1.2pt\def\arraystretch{1.4}
 {\bf J}_{q_0,k} = 
\left[ \begin{array}{c}  
J_{q_0,k,\db_{1}} \\
J_{q_0,k,\db_{2}} \\
\dots \\
J_{q_0,k,\db_{M}}
\end{array} \right].
\label{frechetdef}
\end{equation}
Each block $J_{q_0,k,\db_{m}}$ is the Fr\'{e}chet derivative of the corresponding
mapping $\mathcal{F}_{k,\db_{m}}[q_0]$, whose evaluation in terms of a scattering
problem is described in Theorem \ref{frechet_thm}.
Eq. \eqref{linsingfreqsys} is an overdetermined
linear system of equations for the increment $\delta q$, assuming that 
$M\cdot P$ exceeds the number of degrees of freedom in the representation for
$q(\xb)$, where $P$ denotes the number of equispaced measurements made in the angular
variable on $\parB$. Since we will solve this system iteratively,
we will need an algorithm for applying ${\bf J}_{q_0,k}$ to a vector, as well as its
adjoint ${\bf J}^\ast_{q_0,k}$.

\vspace{.1in}

\begin{theorem} \label{frechet_thm} {\rm \cite{Colton}}
Let $\db$ denote the angle of incidence of an incoming field $u^{inc}$ and let
$u_0 = u^{inc} + u_0^{scat}$ denote the solution to the scattering problem
\begin{equation} 
\Delta u_0(\xb) +k^2(1-q_0(\xb)) u_0(\xb) = 0
\label{u0def}
\end{equation}
in $\mathbb{R}^2$, where $u_0^{scat}$ satisfies the Sommerfeld radiation condition.
Let $\delta q$ be a given perturbation of $q_0$ and let 
$\mathcal{F}_{k,\db}[q_0]$ denote the far field 
operator \eqref{farfielddef}.
Then 
\begin{equation} 
J_{q_0,k,\db}\,  \delta q = v^{far}
\label{Jdef}
\end{equation}
where $v^{far}(\theta) = v(R \cos \theta, R \sin \theta)$ and $v(\xb)$
denotes the solution to the scattering problem
\begin{equation}
\Delta v(\xb) +k^2(1-q_0(\xb)) v(\xb) = k^2 \, \delta q \, u_0 
\label{vdef}
\end{equation}
satisfying the Sommerfeld radiation condition.
\end{theorem}

\vspace{.1in}

\begin{proof}
Let us write the solution to the scattering problem for the inhomogeneity
$q_0 + \delta q$ in the form
\[ 
\Delta (u_0 + v) +k^2(1-q_0 - \delta q)(u_0+v) = 0.
\]
In that case, $v(\xb)$ is the change in the scattered field induced by the 
perturbation $\delta q$.  
The desired result follows after dropping quadratic terms.
\end{proof}

\vspace{.1in}

\begin{theorem} \label{frechet_adj}
Let $f(\theta)$ denote a smooth function on the circle $\parB$ of radius $R$
and let $\chi(f,{\parB})$ denote the corresponding singular charge distribution
on $\parB$ with charge density $f$, viewed as a generalized 
function in $\mathbb{R}^2$.
Let $\db$ denote the direction of incidence of an incoming field $u^{inc}$, and
let $q_0(\xb)$ denote a given inhomogeneity in $\Omega$.
Then the adjoint operator
$J^\ast_{q_0,k,\db}: L^2(\partial \mathcal{B}) \rightarrow C_{0}(\Omega)$ is
given by 
\begin{equation} 
J^\ast_{q_0,k,\db} \, f  = \overline{u_0} \, w
\label{Jastdef}
\end{equation}
where $u_0(\xb)$ denotes the solution to \eqref{u0def} and $w(\xb)$ is the solution to
\begin{equation}
\Delta w(\xb) +k^2(1-q_0(\xb)) w(\xb) =
k^2 \chi(f,{\parB})
\label{wdef}
\end{equation}
in $\mathbb{R}^2$, 
satisfying the adjoint Sommerfeld radiation condition
\[
\lim_{r\rightarrow \infty} \sqrt{r}
\left( \frac{\partial w}{\partial r} +ikw\right) = 0.
\]
\end{theorem}

\vspace{.1in}

\begin{proof}
We first integrate both sides of \eqref{vdef} against the conjugate of $w(\xb)$:
\[ 
\iint_{\mathbb{R}^2}
[\Delta v +k^2(1-q_0)v] \overline{w} dA = 
\iint_\Omega 
k^2 \, \delta q \, u_0 \, \overline{w} dA. 
\]
Using Green's second identity and the Sommerfeld radiation condition,
it is straightforward to show that
\[ 
\iint_{\mathbb{R}^2}
[\Delta \overline{w} +k^2(1-q_0)\overline{w}] v dA = 
\iint_\Omega 
k^2 \, \delta q \, u_0 \, \overline{w} dA. 
\]
or
\[ 
\iint_{\mathbb{R}^2}
 \overline{\chi(f,{\parB})} \, v \, dA = 
\int_\parB 
 \overline{f} \, v \, ds = 
\iint_\Omega 
 \delta q \, u_0 \, \overline{w} dA. 
\]
Since $\langle f,v,\rangle = \int_\parB  \overline{f} \, v \, ds =  
\int_\parB  \overline{f} \, J_{q_0,k,\db} \delta q  \, ds 
= \langle f, J_{q_0,k,\db} \delta q \rangle$, it follows that
\[
\langle J_{q_0,k,\db}^\ast \, f, \delta q \rangle = 
\iint_\Omega 
\, \delta q \, u_0 \, \overline{w} dA =
\langle \overline{u_0} w, \delta q \rangle.
\]
Since $\delta q$ is arbitrary, this yields the desired result.
\end{proof}

\vspace{.1in}

\begin{definition}
We define the adjoint of 
${\bf J}_{q,k}$ by 
\begin{equation}
\arraycolsep=1.2pt\def\arraystretch{1.4}
{\bf J}^\ast_{q_0,k} = 
\left[ 
J^\ast_{q_0,k,\db_{1}} ,
J^\ast_{q_0,k,\db_{2}} ,
\dots ,
J^\ast_{q_0,k,\db_{M}}
 \right].
\label{adjfrechetdef}
\end{equation}
\end{definition}

\subsection{Discretization and regularization at a fixed frequency}

As noted in the introduction, at a given frequency $k$, we can only make $O(k^2)$
independent measurements at finite precision, Thus, we 
seek to reconstruct a model for $q(\xb)$ with 
${\rm supp} (q) \in  [-\frac{\pi}{2},\frac{\pi}{2}]^2$
which has only $O(k^2)$ free parameters.

This avoids various {\em ad hoc} regularization methods that are in common use.
More precisely, at frequency $k$, we approximate the contrast function $q(x_1,x_2)$ 
restricted to the domain $\Omega = [-\frac{\pi}{2},\frac{\pi}{2}]^2$
by the function
\begin{equation}
q_k(x_1,x_2)=
\sum_{\substack{m_1,m_2=1\\ m_1+m_2 \leq S(k)}}^{S(k)} q_{m_1,m_2} 
\sin\left( m_1 \left(x_1 + \frac{\pi}{2} \right) \right) 
\sin\left( m_2 \left(x_2 + \frac{\pi}{2} \right) \right) ,
\label{sineseries}
\end{equation}
with the maximum frequency $S(k) = \lfloor 2k \rfloor$.
This representation has several useful features.
Projection from a sampled function
$q(\xb)$ onto the coefficients $\{ q_{m_1,m_2} \}$ can be accomplished in
$O(N \log N)$ time using the nonuniform FFT (see \cite{DR93,GL2004}
and the references therein). $N$ here denotes the
number of points in the discretization of $q(x_1,x_2)$.
Moreover, the approximation is spectrally accurate for any smooth function $q(\xb)$ 
which has vanished together with all its derivatives at the boundary of $\Omega$. 

\begin{definition}
Let $\hat{q}(k)$ denote the vector of coefficients of the 
truncated sine series in \eqref{sineseries}. We denote by 
$\mathcal{E}_{k}$ the operator which evaluates the 
sine series given by the coefficients $\hat{q}(k)$ at
points $\xb \in \Omega$.
We denote by $\mathcal{E}^\ast_{k}$ its adjoint.
\end{definition}

\subsection{Newton iteration} \label{sec:newton}

Suppose that we have an initial guess $q_k^{(0)}$ for the unknown contrast function,
with far field measurements made at a fixed frequency $k$.
Let $\widehat{\delta q}$ denote the vector of sine series coefficients which we
will use to approximate the unknown perturbation $\delta q$.
Newton's method, for a tolerance $\epsilon$, proceeds as follows:

\vspace{.1in}

\begin{quote}
For $i = 0,1,\dots$
\begin{enumerate}
\item
Solve the linearized problem in a least squares sense using the normal equations:
\begin{equation}
{\mathcal{E}}_{k}^{\ast} 
{\bf J}_{q_k^{(i)},k}^\ast 
{\bf J}_{q_k^{(i)},k} \, 
{\mathcal{E}}_{k}
\widehat{\delta q} =
{\mathcal{E}}_{k}^{\ast} 
 {\bf J}_{q_k^{(i)},k}^\ast  \left(
{\bf u^{\text{\emph{far}}}} - {\bf F}_{k}[q_k^{(i)}] \right).
\label{normaleq}
\end{equation}
\item
Set $q_k^{(i+1)}=q_k^{(i)}+ {\mathcal{E}}_{k} \widehat{\delta q}$.
\item
Stop when 
$\| {\bf u^{\text{\emph{far}}}} - {\bf F}_{k}[q_k^{(i)}] \| < \epsilon$.
\end{enumerate}
\end{quote}

It is instructive, at this stage, to compute the work required at a single 
frequency $k$.  At the $i$th Newton step, we must 
solve $M$ inhomogeneous Helmholtz equations to obtain the right-hand side
for the system \eqref{normaleq}.
Assuming that we solve the normal equations
iteratively using, say, the conjugate gradient method, we must 
solve $2M$ inhomogeneous Helmholtz equations at each iteration
to apply 
${\bf J}_{q_k^{(i)},k}^\ast$ and
${\bf J}_{q_k^{(i)},k}$. 
Each of the PDEs, however, corresponds to a
different right-hand-side in eqs. \eqref{vdef} or \eqref{wdef}. 
Thus, using the HPS solver, we need only compute the factorization of the 
PDE once per Newton iteration.  Thus, the total work is of the order
$O( N_{newton} N^{3/2}) + O(N_{newton} \,(2 N_{iter} +1) \, M \, N)$,
where $N$ denotes the number of grid points used in the solver. 

\vspace{.1in}

\subsection{Recursive Linearization}

Our approach to the full multi-frequency inverse scattering problem 
\eqref{invscatmf} is now straightforward to describe. 
As noted above, it is based on Chen's method of recursive linearization
\cite{Bao_li_2005_1,Bao_li_2005,Bao_li_2009,Bao,Chenrep1088,Chenrep1091,Chen}.

The essential insight of recursive linearization is the following; while
\eqref{invscatmf} is a non-convex, nonlinear system of equations, if a
band-limited approximation $q_k(\xb)$ of $q(\xb)$ were available and $\delta k$
is sufficiently small, then $q_k(\xb)$ is in the basin of attraction for Newton's
method in seeking the global minimum for $q_{k+ \delta k}(\xb)$.
We refer to the references cited above for a discussion of the theoretical 
foundations.
Here, we describe an efficient implementation using all
of the data corresponding to \eqref{invscatmf}.

{\footnotesize
\begin{figure}[htb]
\centering
{\linespread{1.5}
\framebox{\begin{minipage}{14cm} 
{\bf Recursive Linearization using Newton's method}

\vspace{.1in}

We assume we have full aperture data for each of the frequencies
$\{ k_1,k_2,\dots,k_Q\}$ with $k_1 < k_2 < \dots < k_Q$.

\vspace{.1in}

\begin{itemize}
\item Obtain an approximation $q_{k_1}$ for the contrast function $q(\xb)$ at the 
lowest available frequency using the Born approximation
\cite{Bao_li_2005,Bao_Li_2007_NF,Biros1} or a direct imaging method like 
MUSIC or linear sampling \cite{Bao_li_2007}. 
\item For $j=2,\ldots,Q$
\begin{itemize}
\item Create a uniform grid with $N = N(k_j)$ points in the domain $\Omega$. 
\item 
(Since the domain is $k_j/2$ wavelengths across, 10 points per wavelength requires
a grid with $N \approx (5 k_j)^2$ points.)
\item Sample $q_{k_j}$ on the given grid.
\item Solve the single frequency system
${\bf F}_{k_j} [q] = {\bf u}^{far}_{k_j}$ 
using Newton's method \\ (section \ref{sec:newton}) with initial guess
$q_{k_{j-1}}$.
\item Set $q_{k_j}$ to be the solution obtained by Newton's method.
\end{itemize}
\end{itemize}
\end{minipage}}}

\end{figure}
}

A crude estimate of the total work follows, assuming that $k_Q$ is the maximum 
frequency, that we take a step in frequency 
of $\delta k = O(1)$, that the number of Newton iterations $N_{newton} = 1$ and that
the number of iterations $N_{iter}$ required to solve the linear least squares
problem is independent of frequency (see the next section). It is easy to see,
under these hypotheses, that
\[ {\rm Work}\ \approx\ O(k_Q^4) + (2 N_{iter}+1) O(k_Q^4) \, .   \]
The first term is the work required to factor the linear system corresponding to
the forward scattering problem for the initial guess $q_k$ at each successive frequency.
The second term is the work required to solve all the scattering problems required in
applying ${\bf J}_{q_k,k}^\ast$ and ${\bf J}_{q_k,k}$ at each iteration of the 
linearized problem.

\section{Numerical experiments} \label{sec:numerics}

In order to illustrate the performance
of our method, we have chosen four examples of increasing
complexity.
In each case, we take a known function $q(\xb)$ and 
simulate the measured data on $\parB$ by
solving the forward scattering problem. In order to avoid 
``inverse crimes", we use a different solver for data generation than we do 
for inversion. In particular, instead of the HPS solver, we use the fast 
HODLR-based scheme \cite{BorgesLS} for the Lippmann-Schwinger integral equation 
with eight digits of accuracy.

We compute the data 
\[ u^{far}(\theta) = u^{scat}_{q,k_j,\db_{j,m}}(R \cos \theta_p, R \sin \theta_p) \]
for $m=1,\dots,M_j$ at $R=20$ with
$\theta_p=2\pi p/P_j$, for  frequencies 
$k_j=1+j/4$, with $j=0,\ldots,Q$, 
where $M_j=\lfloor2k_j\rfloor$, and
$P_j=\lfloor4k_j\rfloor$.
The incident directions are chosen as 
$\db_{j,m}=\left(\cos \theta_{j,m},\sin \theta_{j,m} \right)$, 
where $\theta_{j,m}=2\pi m/M_j$.

In the first two examples, we use the scattered data computed from our 
forward solver. For the last two examples, 
noise in the form
\begin{equation*}
\tilde{u}^{scat}_{q,k_j,\db{j,m}}(\theta) =
u^{scat}_{q,k_j,\db{j,m}}(\theta) +
\delta \dfrac{\| u^{scat}_{q,k_j,\db{j,m}}(\theta) \|}
{\|\epsilon_1 + i \epsilon_2\|}
(\epsilon_1 + i \epsilon_2)
\end{equation*}
is added, where $\epsilon_1$ and $\epsilon_2$ are normally distributed random variables 
with mean zero and variance one.

For each frequency $k_j$, we discretize the domain $\Omega$ with a uniform quad
tree consisting of $2^l \times 2^l$ square leaf nodes with a
$16 \times 16$ grid used on each to represent $u(\xb)$ and $q(\xb)$.
In examples 1 and 2,
$l$ is chosen so that there are at least 10 points per wavelength in the
discretization, yielding at least 5 digits of accuracy in the solver.
In examples 3 and 4, 
$l$ is chosen so that there are at least 6 points per wavelength in the
discretization, yielding at least 3 digits of accuracy.

For the sake of simplicity, rather than using the 
Born approximation or a direct imaging method
\cite{Biros1,Bao_li_2007,Bao_Li_2007_NF,Bao_li_2005}, we 
assume 
\[
q_1(x_1,x_2)=
\sum_{\substack{m_1,m_2=1\\ m_1+m_2 \leq  2}}^{2} q_{m_1,m_2} 
\sin\left( m_1 \left(x_1 + \frac{\pi}{2} \right) \right) 
\sin\left( m_2 \left(x_2 + \frac{\pi}{2} \right) \right) ,
\]
with $q_{1,1}, q_{1,2}, q_{2,1}$ given as the
projection of $q(\xb)$ onto those modes.

For examples 1--4, we let $k_Q=14.25, 9, 70$, and $70$, respectively.
Finally, we make use of the least squares solver LSQR
\cite{paige1982lsqr} in MATLAB. It is algebraically identical to 
conjugate gradient on the normal equations and the performance 
of the two methods is very similar.
All timings below are reported using 
our solver in conjunction with the parallel computing toolbox in MATLAB,
which makes use of up to 32 cores of a 2.5GHz Intel Xeon system.
Parallelization is straightforward, since the forward scattering problems are all
uncoupled and dominate the CPU time.

\begin{figure}[ht]
  \center
\includegraphics[width=0.65\textwidth]{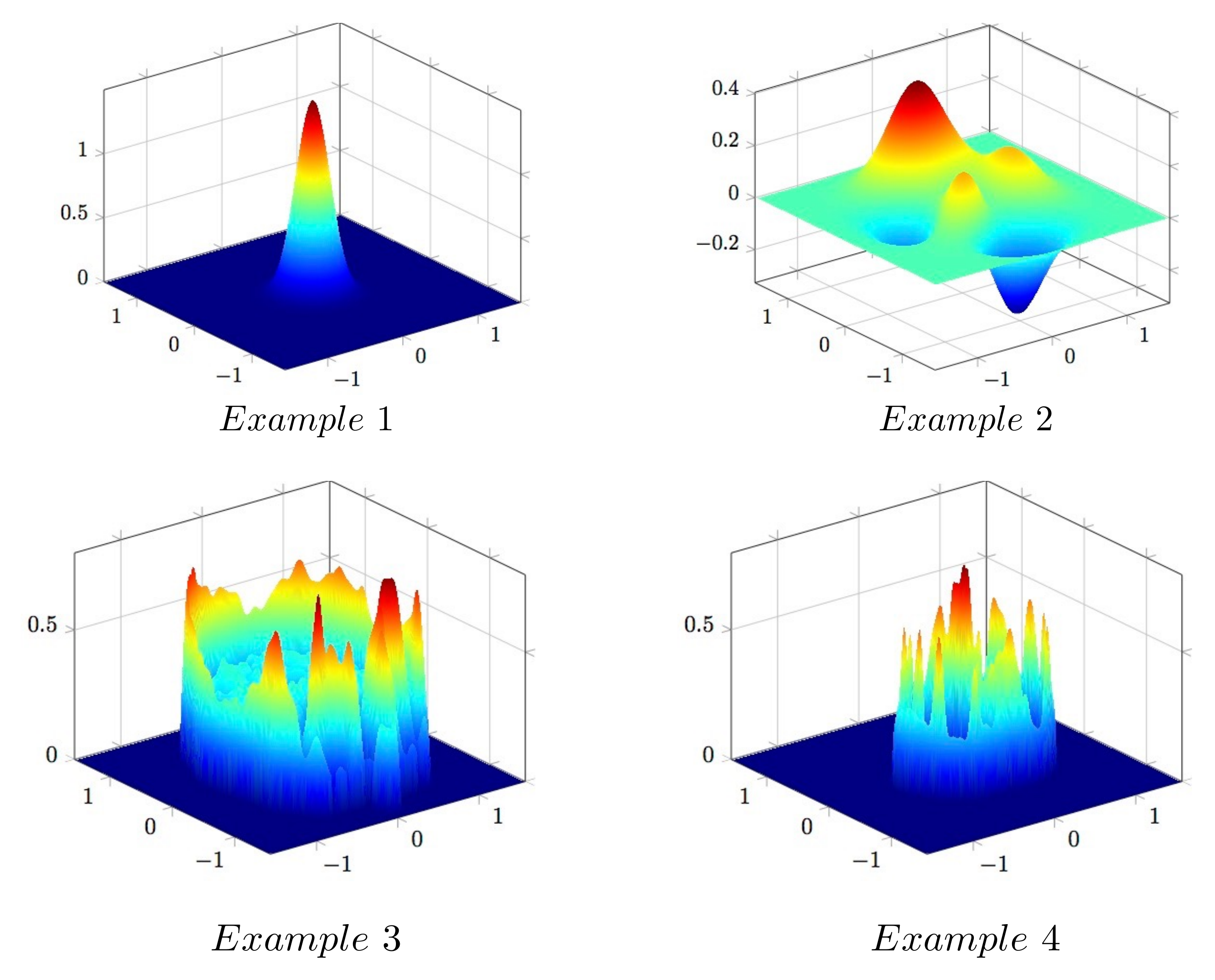}
\caption{The contrast functions for our four examples.}\label{examples_iso}
\end{figure}

\noindent
Example 1: {\em A single Gaussian}.

\begin{figure}[ht]
\center
\includegraphics[width=5in]{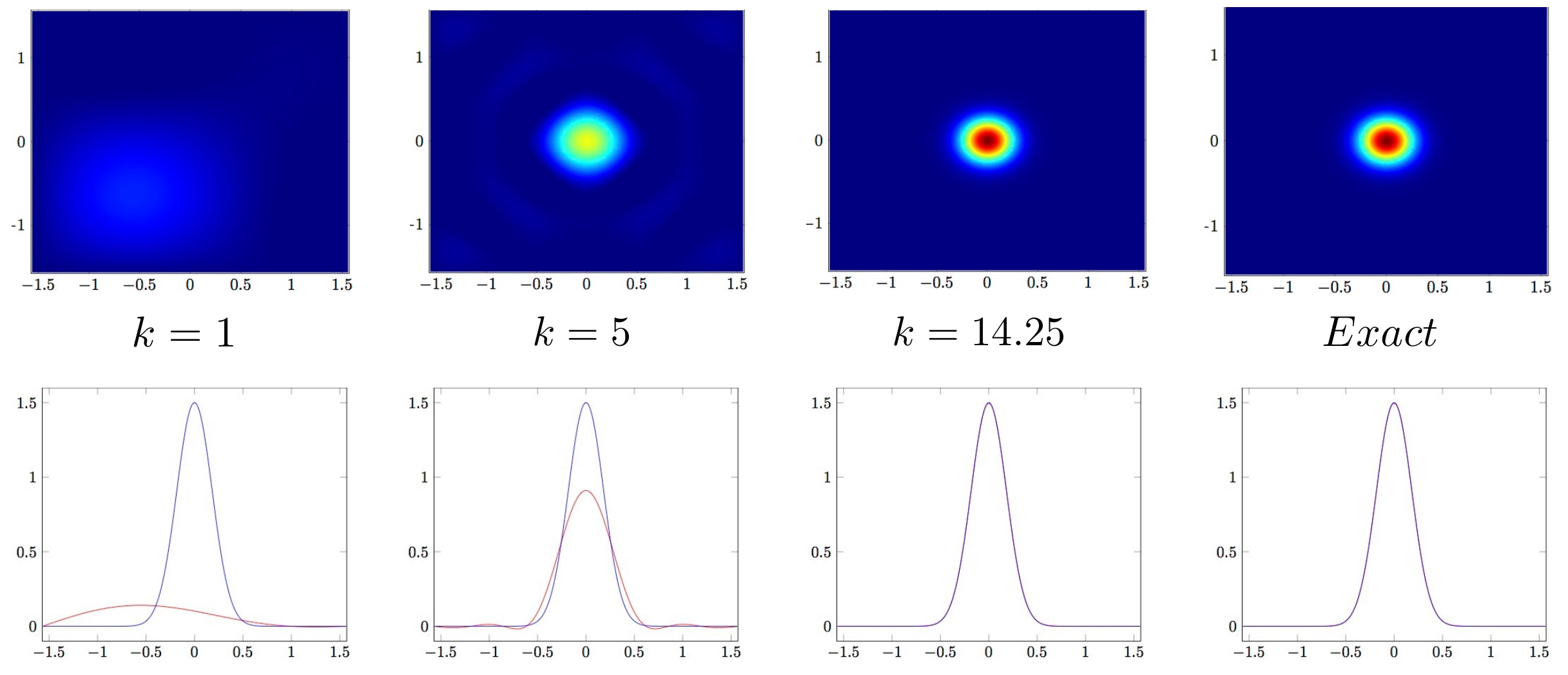}
\caption{\footnotesize Recovering a simple Gaussian contrast function
by recursive linearization (example 1). The upper row shows a contour
plot of the estimated $q(\xb)$ at various frequencies as well as the exact
contrast function.
The lower row shows the corresponding plots of the reconstructed
cross-section $q(x,0)$ for $x\in\left[-\pi/2,\pi/2\right]$.
The reconstruction is shown in red and the original contrast 
is shown in blue.} \label{gaussfig}
\end{figure}

\begin{figure}[ht]
  \center
\includegraphics[width=5in]{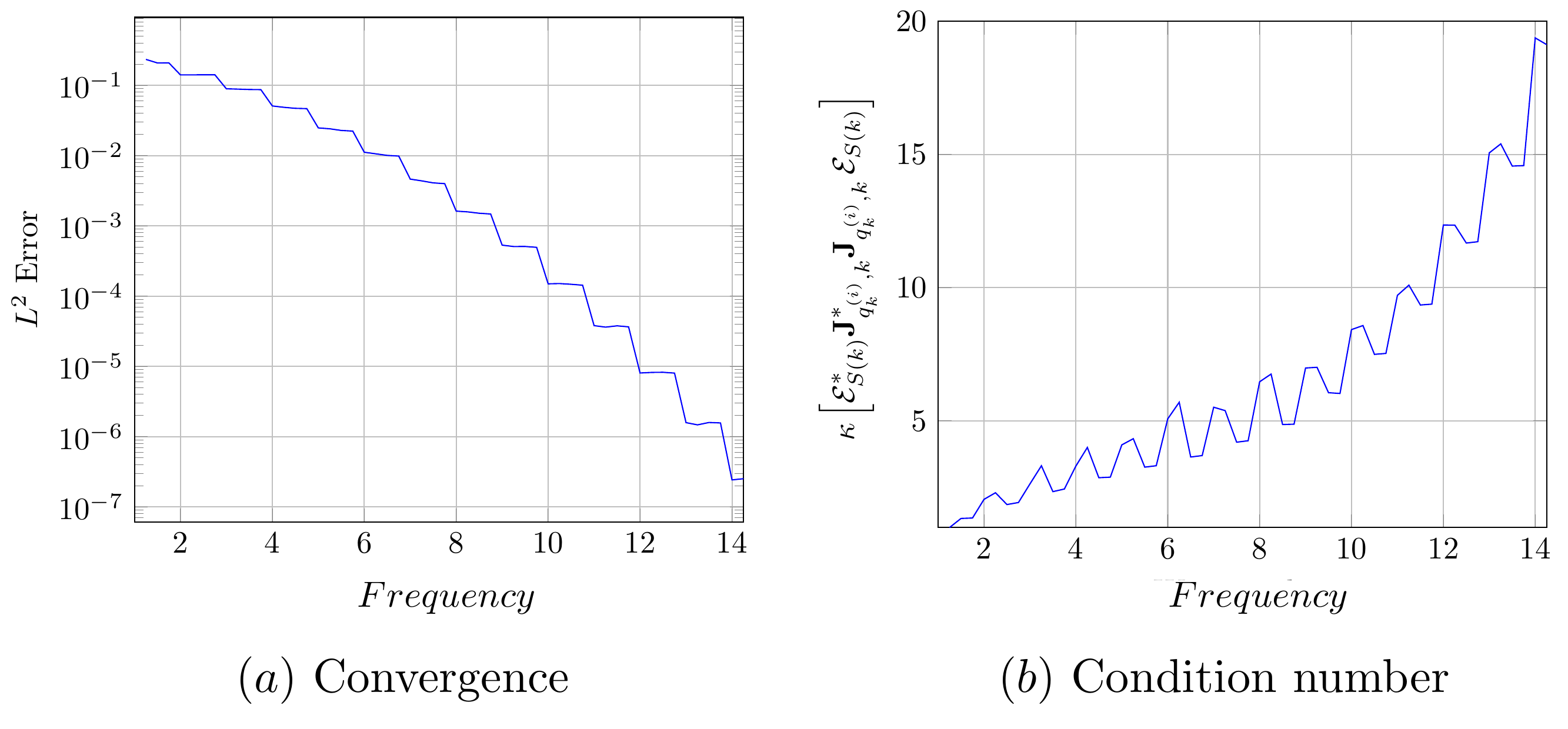}
  \caption{\footnotesize On the left, we plot the error
$\|q_k-q\|$ in $L^2(\Omega)$ as a function of frequency.
On the right, we plot the condition number of the linear least squares
problem as a function of frequency.}       \label{gauss_error}
\end{figure}

First consider the case where the contrast function is a single
Gaussian (Fig. \ref{examples_iso}):
\begin{equation*}
q(x,y)=1.5\exp\left(-\frac{x^2+y^2}{50}\right).
\end{equation*}
The progress of recursive linearization is presented
in Fig. \ref{gaussfig}, which shows contour plots of the exact solution 
next to the reconstructions at the lowest $k = 1$, a mid-range $k = 5$, 
and the highest $k = 14.25$ frequencies.  Below the contour plots are 
cross-sections of the reconstructed function along a single line: that 
is, $q(x,0)$ for $x\in\left[-\pi/2,\pi/2\right]$.  Fig. \ref{gauss_error} 
reports the $L^2$-error of the reconstruction and the condition number of 
the linearized least squares problem as the frequency increases.
 Note that the convergence is very rapid as a function of $k$, since 
 the contrast is smooth and the component solvers are high order accurate.
The total solution time required was about fifteen minutes. 

\noindent
Example 2: {\em A sum of Hermite functions}.

\begin{figure}[ht]
  \center
\includegraphics[width=0.75\textwidth]{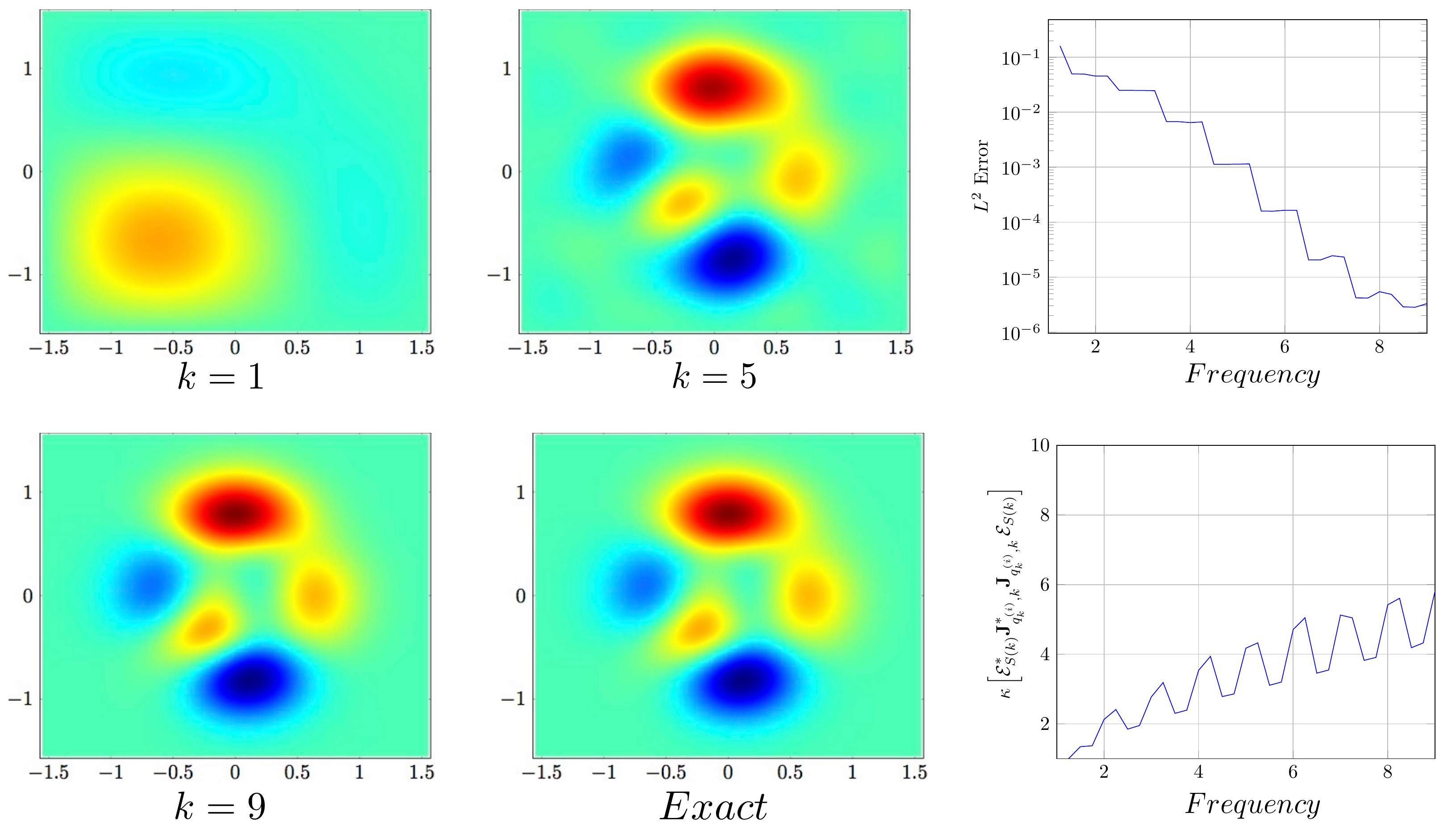}
\caption{\footnotesize 
Recovering a sum of Gaussians
by recursive linearization (example 2). Contour
plots of the estimated $q(\xb)$ are shown at frequencies $k=1$, $k=5$, and $k=9$,
as well as the exact solution.
On the right, we also plot the error
$\|q_k-q\|$ in $L^2(\Omega)$ and the 
the condition number of the linear least squares
problem as functions of frequency.
} \label{fig-ex2}
\end{figure}

We next consider a contrast function made up of a sum of Hermite
functions (Gaussians and their derivatives):
{\footnotesize
\begin{eqnarray*}
q(x,y)&=&0.15\left(1-\frac{x}{\sigma}\right)^2\exp\left(-\left(\left(\frac{x}{\sigma}\right)^2+\left(\frac{y}{\sigma}+1\right)^2\right)\right) -\frac{1}{60}\exp\left(-\left(\left(\frac{y}{\sigma}\right)^2+\left(\frac{x}{\sigma}+1\right)^2\right)\right) \\
&&-\sigma\left(0.4x-\left(\frac{x}{\sigma}\right)^3-\left(\frac{y}{\sigma}\right)^5\right)\exp\left(-\left(\frac{x^2+y^2}{\sigma^2}\right)\right).
\end{eqnarray*}
}
where $\sigma=0.5$ (Fig. \ref{examples_iso}).
While the contrast function in this example is, in some sense,
more complicated than a simple Gaussian, it is a smoother function. Thus,
high fidelity is already achieved at $k=9$.
The progress of recursive linearization is presented
in Fig. \ref{fig-ex2}, which shows contour plots of the 
reconstruction at frequencies $k=1,5,$ and $9$, as well as the exact solution.
The figure also reports the $L^2$-error 
of the reconstruction and the condition 
number of the linearized least squares problem verses frequency. 
Again, the convergence is 
very rapid as a function of $k$, since the contrast is smooth and the 
component solvers are high order accurate.
The total solution time required was about ten minutes.

\noindent
Example 3: {\em Axial cross-section of head}.

\begin{figure}[ht]
\center
\includegraphics[height=6in]{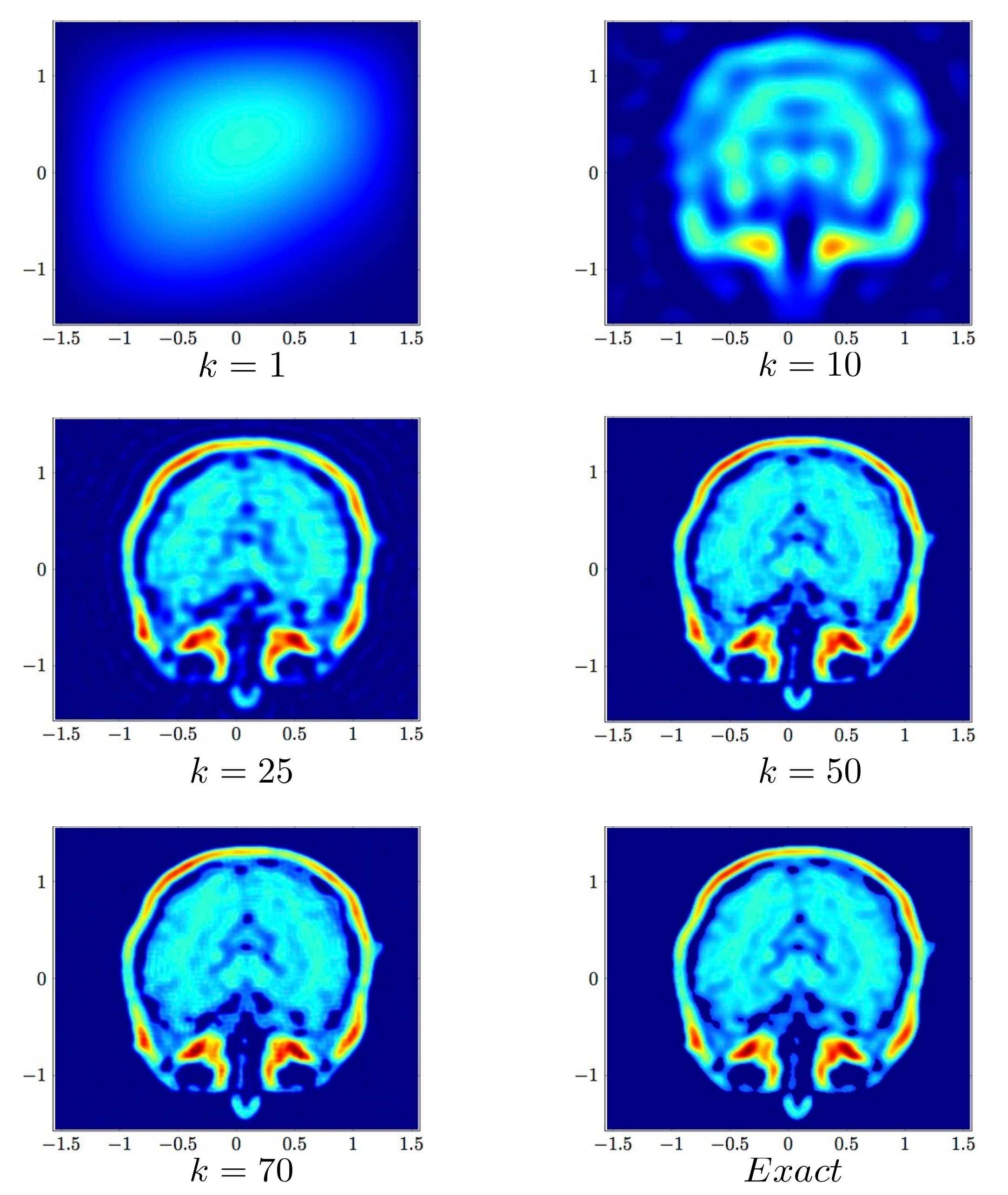} \label{fig-brain}               
\caption{\footnotesize 
Recovering a simulated head phantom 
by recursive linearization (example 3). The 
estimated contrast function $q(\xb)$ is shown at 
frequencies $k=1$, $k=10$, $k=25$, $k=50$, and $k=70$. 
} \label{brain_RLM}       
\end{figure}

For a more interesting (and higher frequency) model, we constructed a
contrast function that resembles the axial
cross section of a human head at the level of the orbitals (a simulated
head phantom). \footnote{The discretized phantom is available 
from the authors upon request.}
A surface plot of the contrast function is shown in Fig. 
\ref{examples_iso} (labeled example 3)
and a contour plot in Fig. \ref{brain_RLM} (labeled Exact). 
Fig. \ref{brain_RLM} also illustrates the
progress of recursive linearization 
at frequencies $k=1, 10, 25, 50,$ and $70$. 
As mentioned previously, our simulated data was computed with 6 
points per wavelength in the discretization, and 5\% noise was added before
reconstruction.

\begin{figure}[ht]
  \center
\includegraphics[height=3.5in]{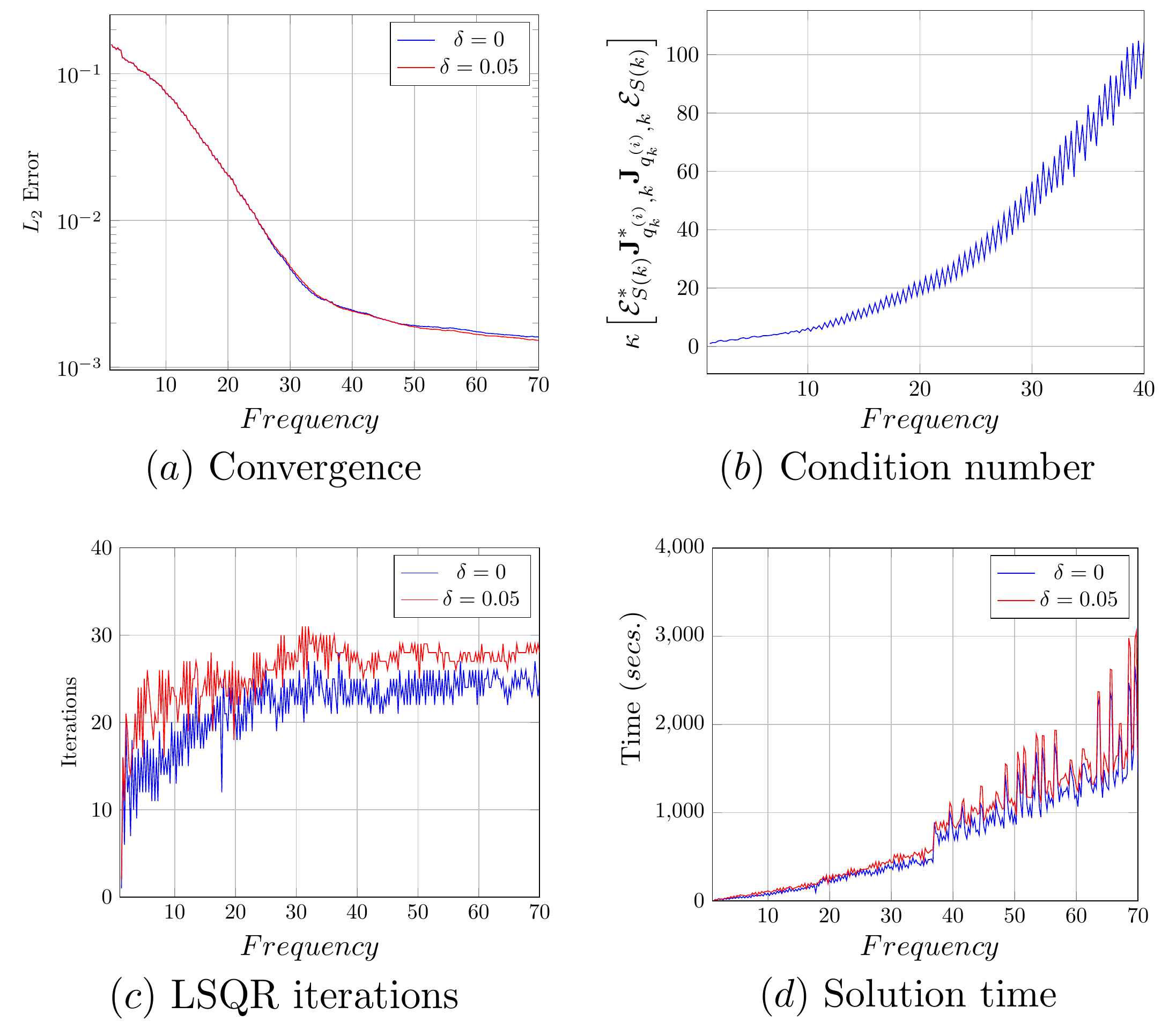} 
  \caption{Numerical results for the simulated axial cross section of a head
(example 3). In (a), we plot the $L^2$ error $\|q_k-q\|$ as a function of frequency.
  In (b), we plot the condition number of the linearized least squares problem and
  in (c) we show the number of LSQR iterations required.
  In (d), we plot the CPU time required at each frequency during the reconstruction
  process.} \label{brainerror} 
\end{figure}

Fig. \ref{brainerror} reports (a) the $L^2$-error 
of the reconstruction, (b) the condition number of the linearized
system, (c) the number of the LSQR iterations required and (d) 
the time in seconds it takes the procedure to create the approximate 
contrast function versus the frequency.  The $L^2$-error, the number of LSQR iterations
and the solution time results are reported for the problem with and without 
$\delta = 0.05$ noise.

Table \ref{tablebrain} provides a more detailed breakdown
of the run time for the recursive procedure (using simulated data
with 5\% noise). Here,
$N$ represents the total number of points used to 
discretize the domain $\Omega$, $Modes = S(k)*[S(k)+1]/2$ is the number of modes 
used as unknowns in the linear least squares problem,
$M$ is the number of incidence directions used, 
$M \cdot P$ is the number of incident directions times the number 
of receiver locations for each $k$,  $T_f$
is the time (in secs.) spent factoring the discretized
forward problem for a given contrast function 
$q_k$, $N_{it}$ is the number of iterations necessary for the LSQR 
method to converge with a tolerance of $\epsilon=10^{-3}$, and $T_l$ 
is the time (in secs.) to solve eq. \eqref{normaleq} at the indicated 
frequency, and $T_t$ is the {\em cumulative} time needed for the full recursion
up to the indicated value of $k$, with steps of $\delta k = 0.25$.

\begin{table}
\begin{center}
\begin{tabular}{c|c|c|c|c|c|c|c|c}
\hline
$k$&
\hspace{0.05cm}$N$\hspace{0.05cm} &
\hspace{0.05cm}$Modes$\hspace{0.05cm} &
\hspace{0.05cm}$M$\hspace{0.05cm} &
$M\cdot P$ &
\hspace{0.05cm}$T_f$\hspace{0.05cm} &
\hspace{0.05cm}$N_{it}$\hspace{0.05cm} &
\hspace{0.05cm}$T_l$ \hspace{0.05cm}& \hspace{0.05cm}$T_t$\\ 
\hline\hline
1.00 & 3721 & 1 & 2 & 16 &    7.77 &
11 &    6.81 &  12.34  \\ 
2.00 & 3721 & 6 & 4 & 64 &    2.95 &
20 &   20.82 &  59.97 \\ 
4.00 & 3721 & 28 & 8 & 256 &    2.93 &
23 &   44.58 &  233.36 \\ 
8.00 & 3721 & 120 & 16 & 1024 &    2.90 &
25 &   92.41 &  918.20 \\ 
16.00 & 3721 & 496 & 32 & 4096 &    2.99 &
27 &  195.14 &  4109.56 \\ 
32.00 & 14641 & 2016 & 64 & 16384 &    6.44 &
30 &  515.72 &  22603.09 \\ 
64.00 & 58081 & 8128 & 128 & 65536 &   25.97 &
28 & 1412.77 & 151114.44 \\ 
\end{tabular}
\end{center}
\caption{Performance of recursive linearization for the simulated
head phantom.}\label{tablebrain}
\end{table}

\noindent
Example 4: {\em Axial cross-section of thorax}.

\begin{figure}[ht]
  \center
\includegraphics[height=6in]{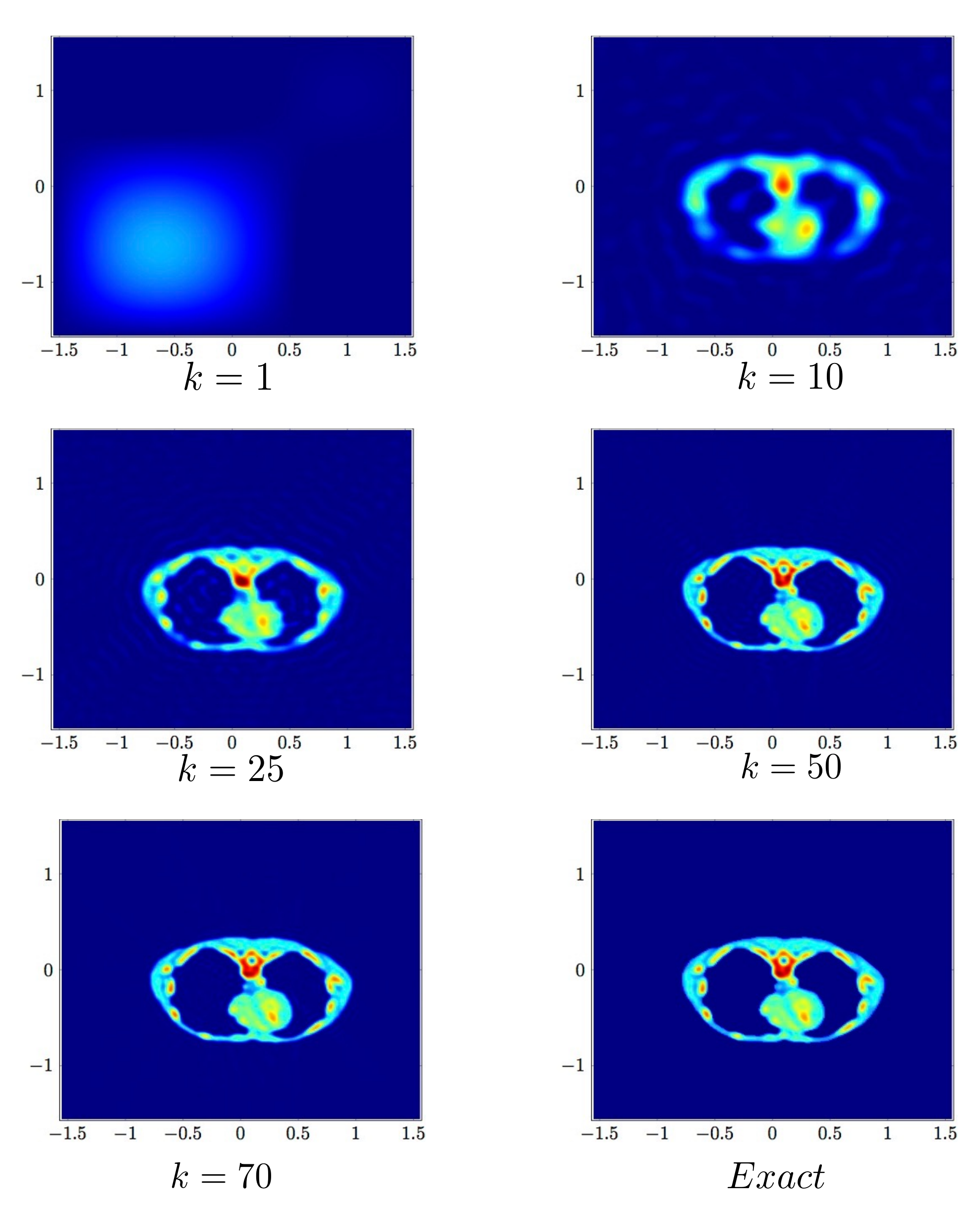} 
  \caption{\footnotesize 
Recovering a simulated thorax phantom 
by recursive linearization (example 4). The 
estimated contrast function $q(\xb)$ is shown at 
frequencies $k=1$, $k=10$, $k=25$, $k=50$, and $k=70$.} \label{fig-thorax}       
\end{figure}

For our last example, we constructed a
contrast function that simulates the axial
cross section of a human thorax at the level of the heart.
\footnote{The discretized phantom is available 
from the authors upon request.}
A surface plot of the contrast function is shown in Fig. \ref{examples_iso} 
(labeled Example 4) and a contour plot
in Fig. \ref{fig-thorax} (labeled Exact). Fig. \ref{fig-thorax} also shows the
progress of recursive linearization 
at frequencies $k=1, 10, 25, 50,$ and $70$. 
The simulated data was computed with 6 
points per wavelength in the discretization, and we added 5\% noise before
reconstruction.

\begin{figure}[ht]
  \center
\includegraphics[height=3.5in]{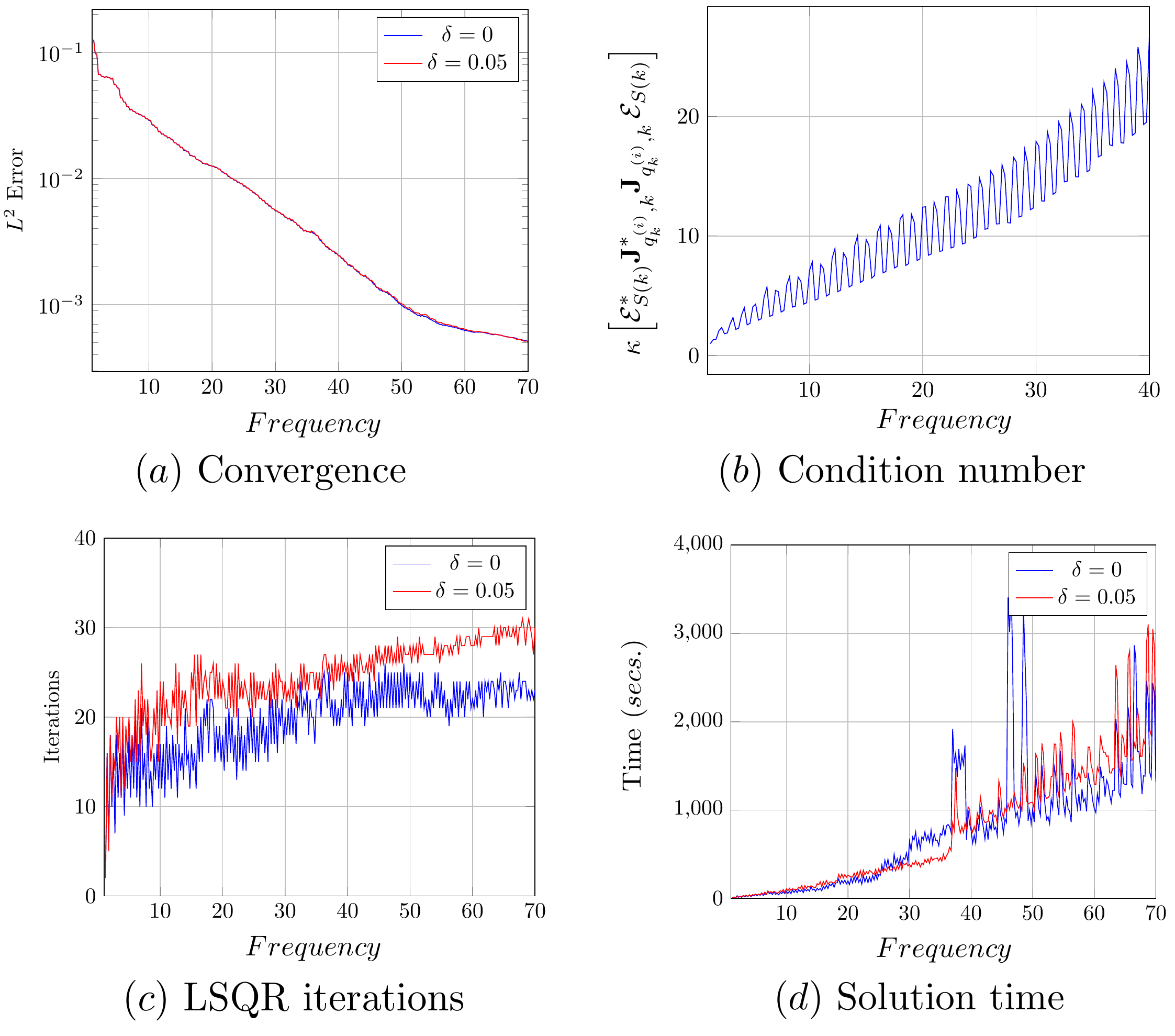} 
\caption{Numerical results for the simulated axial cross section of a thorax
(example 4). In (a), we plot the $L^2$ error $\|q_k-q\|$ as a function of frequency.
  In (b), we plot the condition number of the linearized least squares problem and
  in (c) we show the number of LSQR iterations required.
  In (d), we plot the CPU time required at each frequency during the reconstruction
  process.} \label{thoraxerror} 
\end{figure}

Fig. \ref{thoraxerror} reports (a) the $L^2$-error 
of the reconstruction, (b) the condition number of the linearized
system, (c) the number of the LSQR iterations required and (d) 
the time in seconds it takes the procedure to create the approximate 
contrast function versus the frequency.  The $L^2$-error, the number of LSQR iterations
and the solution time results are reported for the problem with and without 
$\delta = 0.05$ noise.

Table \ref{tablethorax} reports a more detailed breakdown
of the run time for the recursive procedure (using simulated data
with 5\% noise). The same notation is used as in Example 3.

\begin{table}
\begin{center}
\begin{tabular}{c|c|c|c|c|c|c|c|c}
\hline
$k$&
\hspace{0.05cm}$N$\hspace{0.05cm} &
\hspace{0.05cm}$Modes$\hspace{0.05cm} &
\hspace{0.05cm}$M$\hspace{0.05cm} &
$M\cdot P$ &
\hspace{0.05cm}$T_f$\hspace{0.05cm} &
\hspace{0.05cm}$N_{it}$\hspace{0.05cm} &
\hspace{0.05cm}$T_l$\hspace{0.05cm} & $T_t$\\ 
\hline\hline
1.00 & 3721 & 1 & 2 & 16 &    7.57 &
8 &    5.56 &  22.73 \\ 
2.00 & 3721 & 6 & 4 & 64 &    2.73 &
17 &   17.47 & 78.08 \\ 
4.00 & 3721 & 28 & 8 & 256 &    3.12 &
19 &   37.19 &  258.55 \\ 
8.00 & 3721 & 120 & 16 & 1024 &    2.95 &
21 &   77.81 &  995.52 \\ 
16.00 & 3721 & 496 & 32 & 4096 &    3.07 &
26 &  188.05 &  3638.91 \\ 
32.00 & 14641 & 2016 & 64 & 16384 &    6.48 &
23 &  391.33 &  24038.52 \\ 
64.00 & 58081 & 8128 & 128 & 65536 &   25.07 &
29 & 1469.23 & 172365.71 \\ 
\end{tabular}
\end{center}
\caption{Performance of recursive linearization for the simulated
thorax phantom.} \label{tablethorax}
\end{table}

\section{Conclusions} \label{sec:concl}

We have presented a fast, stable algorithm for inverse scattering: reconstructing
an unknown sound speed from far field measurements of the scattered field, in
a fully nonlinear regime.
For this, we have combined Chen's method of recursive linearization with a 
recently developed, spectrally accurate fast direct solver 
\cite{Gillman2014}.
A remarkable feature of
recursive linearization is that by solving a sequence of linearized problems
for sufficiently small steps in frequency (for a commensurate, band-limited
model), one avoids the difficulties associated with the fact that
the high-frequency problem is non-convex and ill-posed.
Using the HPS solver of \cite{Gillman2014}, the CPU time requirements for our scheme
are modest and we believe that the reconstructions shown here are among the largest
ever computed. It is worth noting that for the two large-scale problems
considered above, approximately one million partial differential equations were solved,
requiring approximately two days in our current parallel MATLAB
implementation (using up to 30 cores).

In our experiments, Newton's method requires several iterations at the lowest 
frequency, when the initial guess is far from the desired minimum. 
As the frequency increases, however, a single Newton iteration is sufficient,
consist with the underlying theory \cite{Bao,buithanh2012analysis,Chen}.

Recursive linearization is easily
extended to acoustic or electromagnetic scattering three dimensions. 
All aspects of the scheme described above have clear three-dimensional
analogs. Fast, direct solvers, however, are still under active development
and the scale of the problem is substantially larger, of course, for a fixed 
resolution in each linear dimension.

The scheme described here can be improved and accelerated
in various ways and serves mainly as a ``proof of concept".
Two important issues we have not addressed concern limitations on the available
data; in many settings, only partial aperture data is available and 
in many regimes, only the magnitude of the scattered field can be measured, not its
phase. We are currently working on extensions of the method to such problems.

\section*{Acknowledgments}
This work was supported in part by the 
Applied Mathematical Sciences Program of the U.S. Department of Energy 
under contract DEFGO288ER25053 and by the Office of the Assistant Secretary 
of Defense for Research and Engineering and AFOSR under NSSEFF program award 
FA9550-10-1-0180.
The authors would like to thank Alex Barnett, Yu Chen, Omar Ghattas, Jun Lai, Michael O'Neil, 
Georg Stadler and Tan Bui-Thanh for several useful conversations.
\bibliography{./Bibnew}

\end{document}